\documentclass[12pt, a4paper]{amsart}
\usepackage{amsmath}
\usepackage{geometry,amsthm,graphics,tabularx,amssymb,shapepar}
\usepackage{latexsym,amsfonts,amssymb,amsmath,amsxtra,bbm,color,mathrsfs}
\usepackage{amscd}
\usepackage{hyperref}
\usepackage{pdfsync}
\usepackage[all]{xypic}
\newcommand{\sk}{\mathsf{k}}

\newcommand{\BN}{{\mathbb {N}}}

\newcommand{\CF}{{\mathcal {F}}}

\newcommand{\CS}{{\mathcal {S}}}

\newcommand{\GL}{{\mathrm{GL}}}

\newcommand{\Hom}{{\mathrm{Hom}}}

\newcommand{\Ind}{{\mathrm{Ind}}}

\renewcommand{\Re}{{\mathrm{Re}}}

\newcommand{\rk}{{\mathrm{k}}}

\newcommand{\con}{\textit{C}}

\newcommand{\diag}{\operatorname{diag}}

\newcommand{\sgn}{\operatorname{sgn}}

\newcommand{\od}{\operatorname{d}}

\newcommand{\oL}{\operatorname{L}}

\newcommand{\oZ}{\operatorname{Z}}

\newcommand{\g}{\mathfrak g}

\newcommand{\n}{\mathfrak n}

\renewcommand{\l}{\mathfrak l}

\renewcommand{\rk}{\mathsf k}

\newcommand{\C}{\mathbb{C}}
\newcommand{\R}{\mathbb R}

\newcommand{\abs}[1]{\lvert#1\rvert}
\newcommand{\absk}[1]{\lvert#1\rvert_{\mathsf k}}

\newcommand{\be}{\begin {equation}}
\newcommand{\ee}{\end {equation}}
\newcommand{\bee}{\begin {equation*}}
\newcommand{\eee}{\end {equation*}}

\newcommand{\cf}{\emph{cf.}~}

\theoremstyle{Theorem}

\theoremstyle{Theorem}

\theoremstyle{Theorem}

\theoremstyle{Theorem}
\newtheorem{prp}{Proposition}[section]

\newtheorem{lemp}[prp]{Lemma}

\theoremstyle{Plain}

\theoremstyle{remark}

\theoremstyle{remark}

\theoremstyle{Definition}
\newtheorem{dfn}{Definition}[section]

\newtheorem{prpd}[dfn]{Proposition}
\newtheorem{thmd}[dfn]{Theorem}
\newtheorem{lemd}[dfn]{Lemma}

\newtheorem{remarkd}[dfn]{Remark}

\begin{document}

\title[Rankin-Selberg convolutions]{Rankin-Selberg convolutions for $\GL(n)\times \GL(n)$ and $\GL(n)\times \GL(n-1)$ for principal series representations}

\author[J.-S. Li]{Jian-Shu Li}
\address{Institute for Advanced Study in Mathematics, Zhejiang University, Hangzhou, Zhejiang, P. R. China}
\email{jianshu@zju.edu.cn}

\author[D. Liu]{Dongwen Liu}
\address{School of Mathematical Sciences, Zhejiang University, Hangzhou, Zhejiang, P. R. China}
\email{maliu@zju.edu.cn}

\author[F. Su]{Feng Su}
\address{Department of Pure Mathematics, Xi'an Jiaotong--Liverpool University, Suzhou, Jiangsu, P. R. China}
\email{feng.su@xjtlu.edu.cn}

\author[B. Sun]{Binyong Sun}
\address{Institute for Advanced Study in Mathematics, Zhejiang University, Hangzhou, Zhejiang, P. R. China}
\email{sunbinyong@zju.edu.cn}

%\address{Academy of Mathematics and Systems Science, Chinese Academy of
%Sciences \& University of Chinese Academy of Sciences,  Beijing, 100190, China} \email{sun@math.ac.cn}


%\thanks{J.L.'s research was partially supported by RGC-GRF grant 16303314 of HKSAR}
%\thanks{Supported by NSFC Grant 11222101}

\begin{abstract}
Let $\rk$ be a local field. Let $I_\nu$ and $I_{\nu'}$ be smooth principal series representations of $\GL_n(\rk)$ and $\GL_{n-1}(\rk)$ respectively. The Rankin-Selberg integrals yield a continuous bilinear map $I_\nu\times I_{\nu'}\longrightarrow \C$ with a certain invariance property. We study  integrals over a certain open orbit that also yield a continuous bilinear map $I_\nu\times I_{\nu'}\longrightarrow \C$ with the same invariance property, and show that these integrals equal the Rankin-Selberg integrals up to an explicit constant. Similar results are also obtained for  Rankin-Selberg integrals for $\GL_n(\rk)\times \GL_n(\rk)$.  
%Let $E$ be a quasi-complete Hausdorff locally convex complex topological vector space, and let $M$ be a paracompact smooth manifold. We show that every element of  $\RC^\infty(M;E)$ belongs to the closure of a bounded subset in $\RC^\infty(M)\otimes E\subset  \RC^\infty(M;E)$. Consequently,
% the topological vector space $\RC^\infty(M;E)$ is identical to the quasi-complete hull of the projective tensor product $\RC^\infty(M)\otimes E$.
%Let $\RD$ be a finite dimensional central division algebra over a local field $\RF$. A relation is established between poles of  Godement-Jacquet L-functions for $\GL_n(\RD)$  and  distributions on the matrix space $\oM_n(\RD)$  with  support in  the set of %singular matrices. As an application, we show that the full theta lifts to $\GL_n(\RF)$ of generic irreducible representations of $\GL_n(\RF)$ are irreducible.

%Let $\RD$ be a finite dimensional central division algebra over a local field $\RF$. Let $\sigma$ be an irreducible admissible smooth representation of $\GL_n(\RD)$.
%We use Godement-Jacquet L-functions to give a sufficient condition for the irreducibility of the full theta lift of $\pi$ to $\GL_n(\RD)$. In particular, we show that if $\RD=\RF$ and $\pi$ is generic, then the full theta lift is irreducible.
\end{abstract}

\maketitle

\section{Introduction and the main results}

Although Rankin-Selberg convolution is a well-established theory, explicit calculations of Rankin-Selberg integrals are usually not easy. These explicit calculations are often crucial for the arithmetic study of Rankin-Selberg L-functions.

Let $n$ be a positive integer, and $n':=n$ or $n-1$ throughout this article. The Rankin-Selberg convolutions for $\GL(n)\times \GL(n')$ are viewed as the basic cases of the general  Rankin-Selberg theory.  In these basic cases, at least for principal series representations, we  aim to calculate the Rankin-Selberg integrals as explicitly as possible. More precisely, we will  express the  Rankin-Selberg integral   as a more explicit integral over a certain $\GL(n')$-torsor. 

The archimedean case of the main result (Theorem \ref{thmA}) of this article is used as a key ingredient in \cite{LLS21} to prove the period relations for the critical values of Rankin-Selberg L-functions, which is an automorphic analog of Deligne’s conjecture (\cite{D79}). In a paper under preparation, the non-archimedean case of Theorem \ref{thmA} will be used to calculate the modified Euler factor at place $p$ for $p$-adic Rankin-Selberg L-functions (as predicted  by J. Coates, see \cite{C-PR89,
C89a, C89b}).

\subsection{Principal series representations}
Fix an  arbitrary local field $\sk$. Write $\abs{\,\cdot\,}_\sk: \sk\longrightarrow \R$ for the normalized absolute value. Fix an arbitrary nontrivial unitary character $\psi: \sk\longrightarrow \C^\times$. We equip $\sk$ with the self-dual Haar measure associated to $\psi$.

 For every $k\in \BN:=\{0,1,2,\cdots\}$, write $G_k:=\GL_k(\sk)$. It contains $\bar B_k N_k$ as an open dense subset, where $\bar B_k$ is the subgroup of the lower triangular matrices, and $N_k$ is the subgroup of the unipotent upper triangular matrices.
 The group $N_k$ is equipped with the Haar measure
\be\label{mea1}
\od\! u:=\prod_{1\leq i< j\leq k} \od\!u_{i,j}, \qquad u=[u_{i,j}]_{1\leq i,j\leq k}\in N_k,
 \ee
the group $G_k$ is  equipped with the Haar measure
 \[
 \od\! g:=\abs{\det g}_\sk^{-1} \cdot \prod_{1\leq i, j\leq k} \od\!g_{i,j}, \qquad g=[g_{i,j}]_{1\leq i,j\leq k}\in G_k.
\]
and the  group $\bar B_k$ is equipped with the left invariant Haar measure
\be\label{mea3}
 \od\! \bar b:=\prod_{i=1}^k \abs{\det \bar b_{i,i}}_\sk^{-i} \cdot \prod_{1\leq j\leq i\leq k} \od\! \bar b_{i,j}, \qquad \bar b=[\bar b_{i,j}]_{1\leq i,j\leq k}\in \bar B_k.
\ee
Here $\od\!u_{i,j}$, $\od\!g_{i,j}$ and $\od\! \bar b_{i,j}$ indicate the Haar measure on $\sk$ associate to $\psi$ as before.
The coset $N_k\backslash G_k$ is  equipped with the invariant quotient measure.
Unless otherwise mentioned, all measures appearing in integrals in this article are the specified measures as above.

As usual, a continuous homomorphism from a topological group to $\C^\times$ is called a character of the topological group. Write $\widehat{\sk^\times}$ for the set of all characters of $\sk^\times$.  For every  $\mu\in (\widehat{\sk^\times} )^k$,   we view it as a character of $\bar B_k$ in the obvious way, and write
\[
I_\mu:=\Ind_{\bar B_k}^{G_k}\mu :=\{f\in \con^\infty(G_k)\mid f(\bar b x)=\mu(\bar b)\cdot \bar \rho_k(\bar b)\cdot f(x)\textrm{ for all } \bar b\in \bar B_k, \ x\in G_k\}
 \]
for the corresponding smooth principal series representation, on which $G_k$ acts by right translation. Here
 \[
   \bar \rho_k:=(\abs{\,\cdot\,}_\sk^\frac{1-k}{2}, \abs{\,\cdot\,}_\sk^\frac{3-k}{2}, \cdots, \abs{\,\cdot\,}_\sk^\frac{k-1}{2})\in (\widehat{\sk^\times} )^k.
 \]
  In the archimedean case, $I_\mu$ is naturally a Fr\'echet space.  In the non-archimedean case, $I_\mu$ is countable-dimensional. We view every countable-dimensional complex vector space as a locally convex topological vector space with the finest locally convex topology. In particular, every linear functional on $I_\mu$ is continuous in the  non-archimedean case.

\subsection{Rankin-Selberg integrals}

Write $\CS(X)$ for the space of Schwartz functions on $X$ when $X$ is a Nash manifold (see \cite{AG08}), and the space of compactly supported locally constant functions on $X$ when $X$ is a totally disconnected locally compact topological space. All functions in this article are complex-valued. %unless otherwise specified. 

We review some basic facts concerning the Rankin-Selberg convolutions (see \cite{JacquetPSShalikaRankinSelberg} and \cite{J09} for more details).
 Define a character
\[
  \psi_k: N_k\longrightarrow \C^\times, \quad [u_{i,j}]_{1\leq i,j\leq k}\longmapsto \psi\left(\sum_{1\leq i\leq k-1} u_{i, i+1}\right).
\]
 When no confusion is possible, we will not distinguish a character with the corresponding representation on $\C$.
The space  $\Hom_{N_k}(I_\mu, \psi_k)$ of the $N_k$-equivariant continuous linear functionals is one-dimensional, and there is a unique element  of it, to be denoted by $\lambda_\mu$,   such that (see \cite[Theorem 15.4.1]{Wa92}) 
\[
  \lambda_\mu(f)=\int_{N_k} f(u) \overline{\psi_k}(u)\od\! u
\]
for all $f\in I_\mu$ such that $f|_{N_k}\in \CS(N_k)$. Here and henceforth, an overline over a character indicates its complex conjugation.   Similarly, denote by $\lambda'_\mu$  the unique element of $\Hom_{N_k}(I_\mu, \overline{\psi_k})$ such that
\[
  \lambda'_\mu(f)=\int_{N_k} f(u) \psi_k(u)\od\! u
\]
for all $f\in I_\mu$ such that $f|_{N_k}\in \CS(N_k)$.  Then we have the homomorphisms
\[
 W: I_\mu\longrightarrow \Ind_{N_k}^{G_k} \psi, \quad f\longmapsto W_f:=(g\longmapsto \lambda_\mu(g.f) )
 \]
 and 
\[ \overline{W}: I_\mu\longrightarrow \Ind_{N_k}^{G_k} \overline{\psi}, \quad f\longmapsto \overline{W}_f:=(g\longmapsto \lambda_\mu'(g.f)).
\]

Recall that  $n'=n$ or $n-1$. Throughout this article we fix   
\[
\nu=(\nu_1, \nu_2, \cdots,\nu_n)\in (\widehat{\rk^\times} )^n\qquad\textrm{and}\qquad  \nu'=(\nu'_1, \nu'_2, \cdots,\nu'_{n'})\in (\widehat{\sk^\times} )^{n'}.
\]
  Put
\[
  \oL(s, \nu\times \nu'):=\prod_{1\leq i\leq n, 1\leq j\leq n'} \oL(s, \nu_i\cdot \nu'_j).
\]
Here and as usual, for every character $\chi$ of $\sk^\times$,  $\oL(s,\chi)$ denotes the local L-function of $\chi$.
For all $k,l\in \BN$,  denote by  $\sk^{k\times l}$ the space of $k\times l$ matrices with entries in $\sk$.
Let $f\in I_\nu$, $f'\in I_{\nu'}$ and $\phi \in \CS(\sk^{1\times n})$.

If $n'=n$, the Rankin-Selberg integral is defined by
\be\label{rsnn}
  \oZ(s, f, f', \phi):=\int_{N_n\backslash G_n} W_f(g) \cdot \overline{W}_{f'}(g)\cdot \phi(e_n g)\cdot \abs{\det g}_\sk^s\,\od\! g,
\ee
where $e_n:=[0,0, \cdots, 0,1]\in \sk^{1\times n}$. The integral \eqref{rsnn} is absolutely convergent when the real part $\Re(s)$ of the complex variable $s$ is sufficiently large, and extends to a holomorphic multiple of $\oL(s, \nu\times \nu')$ (see \cite[Section 8.1]{J09}).
More precisely, there exists a unique continuous map
\be\label{zcirc}
  \oZ^\circ: \C\times  I_{\nu}\times I_{\nu'} \times \CS(\sk^{1\times n})\longrightarrow \C
\ee
with the following properties:
\begin{itemize}
\item it is holomorphic in the first variable and linear in the last three variables;
\item there exists a constant $c_{\nu,\nu'}\in \R$ such that whenever $\Re(s)>c_{\nu,\nu'}$, the integral \eqref{rsnn} is  absolutely convergent  and
\[
  \oZ(s, f, f', \phi)=\oL(s,\nu\times \nu')\cdot \oZ^\circ(s,f, f', \phi )
\]
for all $f\in I_\nu$, $f'\in I_{\nu'}$  and $\phi\in \CS(\sk^{1\times n})$.
\end{itemize}

If $n'=n-1$, the Rankin-Selberg integral is defined by
\be\label{rsnnp}
  \oZ(s, f, f'):=\int_{N_{n-1}\backslash G_{n-1}} W_f(g)\cdot \overline{W}_{f'}(g)\cdot  \abs{\det g}_\sk^{s-\frac{1}{2}}\,\od\! g.
\ee Similarly, the integral \eqref{rsnnp} is absolutely convergent when  $\Re(s)$ is sufficiently large, and extends to the multiplication of $\oL(s, \nu\times \nu')$ with a
continuous map
\be\label{zcirc}
  \oZ^\circ: \C\times  I_{\nu}\times I_{\nu'} \longrightarrow \C
\ee
that is holomorphic in the first variable and linear in the last two variables.

\subsection{The integrals over the open orbits}
The right action of $G_{n}$ on $(\bar B_n\backslash G_n)\times (\bar B_{n}\backslash G_{n})\times \rk^{1\times n}$ has a unique open orbit. 
Likewise, the  right action of $G_{n-1}$ on $(\bar B_n\backslash G_n)\times (\bar B_{n-1}\backslash G_{n-1})$ has a unique open orbit.
We will introduce  integrals over these open orbits and relate them to the Rankin-Selberg convolutions. For this purpose, we introduce some auxiliary matrices as follows.
For each $k\in \BN$, write
\[
 w_k:=\left[
        \begin{array}{cccc}
          0 & \cdots & 0 & 1 \\
          0 & \cdots & 1 & 0 \\
          &\cdots & \cdots &  \\
          1 & 0 & \cdots & 0 \\
        \end{array}
      \right]\in  \GL_k(\sk).
\]
Define a family  $\{z_k\in \GL_k(\sk)\}_{k\in \BN}$ of matrices
inductively by
\[
 z_0:=\emptyset\ \  (\textrm{the unique element of $\GL_0(\sk)$}), \quad  z_1:=[1],
 \]
and
\[
       z_k :=
 \left[
            \begin{array}{cc}
                            w_{k-1}& 0 \\
                       0 & 1 \\
                     \end{array}
                   \right]
                    \left[
            \begin{array}{cc}
                       z_{k-2}^{\iota}& 0 \\
                       0 & 1_2 \\
                     \end{array}
                   \right]
                    \left[
            \begin{array}{cc}
                       {}^tz_{k-1}  w_{k-1} z_{k-1}& {}^t e_{k-1}\\
                       0 & 1 \\
                     \end{array}
                   \right], \quad \textrm{for all  $k\geq 2$.}
                   \]
Here and as usual, a left superscript $t$ over a matrix  indicates the transpose, a right superscript $\iota$  indicates the inverse transpose of an invertible matrix,
$1_2$ stands for the $2\times 2$ identity matrix, and $e_{k-1}:=[0,\cdots, 0,1]\in \sk^{1\times(k-1)}$. In particular, 
$z_2:=\left[
  \begin{array}{cc}
    1 & 1  \\
     0& 1 \\
                  \end{array}
\right]$.

\begin{lemd} \label{lem:orbit}
{\rm (a)} The right action of $G_{n}$ on $(\bar B_n\backslash G_n)\times (\bar B_{n}\backslash G_{n})\times \rk^{1\times n}$ has a unique open orbit represented by 
\begin{equation}\label{repr1}
\left(z_n, \begin{bmatrix} z_{n-1} & 0 \\ 0 & 1\end{bmatrix}, e_n\right).
\end{equation}
{\rm (b)} The right action of $G_{n-1}$ on $(\bar B_n\backslash G_n)\times (\bar B_{n-1}\backslash G_{n-1})$ has a unique open orbit represented by
\begin{equation}\label{repr2}
(z_n, z_{n-1}).
\end{equation}
\end{lemd}

\begin{proof} We first prove  inductively the claim that the stabilizer (in $G_n$ or $G_{n-1}$) of the element represented by \eqref{repr1} or \eqref{repr2} is trivial. The claim is trivial for $n=1$, and we assume that $n\geq 2$. Since the stabilizer of $e_n$ in $G_n$ is the mirabolic subgroup 
\[
\left\{\left.\begin{bmatrix} g & v \\ 0 & 1 \end{bmatrix} \right| \, g\in G_{n-1}, \, v\in \rk^{(n-1)\times 1}\right\},
\]
it is easy to see that the claim in case (a) follows from that of case (b). 

We will show that the claim in case (b) follows from the validity of the claim in case (a) for $n-1$, which thereby finishes the proof by induction. Consider the diagonal action of $G_{n-1}$ on 
\[
(\bar B_{n-1} \backslash G_{n-1})\times (\bar B_{n-1}\backslash G_{n-1})\times \rk^{(n-1)\times 1},
\]
where the right action of  $g\in G_{n-1}$ on $\rk^{(n-1)\times 1}$ is given by
$v\mapsto g^{-1}v$. Direct computation shows that the stabilizer $H$ of $(\bar B_n z_n, \bar B_{n-1}z_{n-1})$ in $G_{n-1}$ is contained in the stabilizer of 
\[
\left(\bar B_{n-1} w_{n-1} \begin{bmatrix} 
z_{n-2}^\iota & 0 \\ 0 & 1\end{bmatrix} {}^tz_{n-1} w_{n-1} z_{n-1}, \bar B_{n-1} z_{n-1}, z_{n-1}^{-1}w_{n-1}  z_{n-1}^\iota {}^te_{n-1}\right)
\]
under the above action. The latter stabilizer is conjugate to the stabilizer of 
\[
\left(\bar B_{n-1} w_{n-1} \begin{bmatrix} 
z_{n-2}^\iota & 0 \\ 0 & 1\end{bmatrix}, \bar B_{n-1} w_{n-1}  z_{n-1}^\iota, {}^te_{n-1}\right),
\]
which is the matrix transpose of the stabilizer $H'$ of 
\[
\left( \bar B_{n-1}\begin{bmatrix} 
z_{n-2} & 0 \\ 0 & 1\end{bmatrix},  \bar B_{n-1} z_{n-1}, e_{n-1}\right).
\]
If the claim in case (a) holds for $n-1$, then the group $ H'$ is trivial, hence $H$ is trivial as well, which implies that the claim in case (b) holds for $n$. 

For both (a) and (b), dimension counting shows that an orbit is open if and only if the stabilizers are  finite groups. In particular, the above argument implies that the element represented by \eqref{repr1} or \eqref{repr2} belongs to an open orbit.

We next prove the uniqueness of the open orbit (This is known to experts, and we include a proof for completeness.).

For (a), if a $G_n$-orbit  in $(\bar B_n\backslash G_n)\times (\bar B_{n}\backslash G_{n})\times \rk^{1\times n}$ is open, then its image in $(\bar B_n\backslash G_n)\times (\bar B_{n}\backslash G_{n})$ under the natural projection is also open, which has to be the orbit of $(\bar B_n, \bar B_n w_n)$, as is well-known. The stabilizer of $(\bar B_n,\bar B_n w_n)$ in $G_n$ is the diagonal torus, whose action on $\rk^{1\times n}$ has a unique open orbit  $(\rk^\times)^{1\times n}$. This shows the uniqueness in case (a).

For (b), note that we have a $G_{n-1}$-equivariant open embedding 
\[
\bar B_{n-1}\backslash G_{n-1}\times \rk^{1\times(n-1)}\hookrightarrow \bar B_n \backslash G_n,\quad (\bar B_{n-1} g, v)\mapsto  \bar B_n w_n \begin{bmatrix} w_{n-1}g & 0 \\   v & 1\end{bmatrix}
\]
with dense image. Then the uniqueness of the open orbit in (b) follows from applying the uniqueness assertion in (a) for $n-1$.
\end{proof}

We have some remarks for the elements $z_k$. Note that these elements are rational, which is required in the study of period relations in \cite{LLS21}. Below we will introduce certain integrals over the open orbits in Lemma 
\ref{lem:orbit}. The inductive choice of $z_k$, which looks complicated at first glance, will yield nice  recurrence relations for these integrals in Section \ref{sec3}.

We are concerned with  the following two integrals.

\begin{dfn} Let $f\in I_\nu$, $f'\in I_{\nu'}$ and $\phi \in \CS(\sk^{1\times n})$. 

\noindent
(a) 
Suppose that $n'=n$. For every $s\in \C$, define
\be\label{nn}
 \Lambda(s, f, f', \phi):=\int_{G_{n}} f(z_n g
 )\cdot f'\left(\left[
                                                                \begin{array}{cc}
                                                                  z_{n-1} & 0 \\
                                                                  0 & 1 \\
                                                                \end{array}
                                                              \right]  g\right) \cdot \phi(e_n g)\cdot \abs{\det g}^{s}_{\rk}\od\! g.
\ee

\noindent
(b) 
Suppose that $n'=n-1$. For every $s\in \C$, define
\be\label{nnp}
 \Lambda(s, f, f'):=\int_{G_{n-1}} f\left(z_n \left[
                                                                \begin{array}{cc}
                                                                  h & 0 \\
                                                                  0 & 1 \\
                                                                \end{array}
                                                              \right] \right)\cdot f'(z_{n-1} h) \cdot \abs{\det h}^{s-\frac{1}{2}}_{\rk}\od\! h.
\ee

\end{dfn}

The following lemma is clear.

\begin{lemd} \label{lem1.1} Let  $f\in I_\nu$, $f'\in I_{\nu'}$ and $\phi \in \CS(\sk^{1\times n})$. 

\noindent {\rm (a)}
 If $n'=n$, then for $g\in G_n$, 
\be\label{linv}
\Lambda(s, g. f, g. f', g.\phi )=\abs{\det g}^{-s}_{\rk} \Lambda(s, f, f', \phi)
\ee
and
\[
  \oZ(s, g. f, g. f', g.\phi)=\abs{\det g}^{-s}_{\rk} \oZ(s,f, f', \phi).
\]

\noindent {\rm (b)}
 If $n'=n-1$, then for $h\in G_{n-1}$,
\[
\Lambda(s, h. f, h. f' )=\abs{\det h}^{-s+\frac{1}{2}}_{\rk} \Lambda(s, f, f')
\]
and
\[
  \oZ(s, h. f, h. f')=\abs{\det h}^{-s+\frac{1}{2}}_{\rk} \oZ(s,f, f').
\]

\end{lemd}

More precisely, the left hand side integral of \eqref{linv} is absolutely convergent if and only if so is the right one, and when this is the case the equality \eqref{linv} holds.
Similar interpretation applies to other equalities for integrals or double integrals in this article (for example, in Proposition \ref{r1} and \ref{prprl1}) without further comments.

For every character $\omega$ of $\sk^\times$, denote by $\mathrm{ex}(\omega)$ the unique real number such that
\[
  \abs{\omega(a)}=\absk{a}^{\mathrm{ex}(\omega)},\quad \textrm{ for all }a\in \sk^\times.
\]

\begin{prpd} \label{prop1.2}

Assume that $s$ lies in the vertical strip 
\[
\Omega_{\nu, \nu'}:=\left\{ s\in \mathbb{C}\left|\begin{aligned} & \mathrm{ex}(\nu_i)+\mathrm{ex}(\nu'_j)+\Re(s)<1  \textrm{ whenever } i+j\leq n,\\
& \mathrm{ex}(\nu_i)+\mathrm{ex}(\nu'_j)+\Re(s)>0 \textrm{ whenever } i+j> n
\end{aligned}\right.\right\}.
\]
Then the integrals \eqref{nn} and   \eqref{nnp} are absolutely convergent.
\end{prpd}

We remark that the set $\Omega_{\nu, \nu'}$ in the above proposition may or may not be empty. 

% for all $f\in I_\mu$ whose support has the form $\bar B_n C$, where $C$ is a compact subset of $N_n$.

\subsection{An equality of  two integrals} \label{sec1.4}

%For every $\mu=(\mu_1,\mu_2, \cdots, \mu_k)\in (\widehat{\sk^\times})^k$, $k\in \BN$, write
% \begin{eqnarray*}
  %\sgn(\mu)&:=&\sgn(\mu_1,\mu_2, \cdots, \mu_k):=\prod_{1\leq i \leq \frac{k+1}{2} \textrm{ and $2|i-1$, or $\frac{k+1}{2}\leq i \leq k$ and $2|k-i$}} \mu_i(-1)^k.
%\end{eqnarray*}
%For $\nu\in (\widehat{K^\times})^n$ and $\nu'\in (\widehat{K^\times})^{n'}$, define
%\[
% \sgn(\nu;\nu'):=                     \sgn(\nu_2,\nu_3, \cdots, \nu_{n-1})\cdot \sgn(\nu'_1,\nu'_2, \cdots, \nu'_{n-1}),
%\]
%which by convention equals $1$ when $n=1$.

Define a sign 
\[
 \sgn(\nu;\nu'):=\prod_{j<i, \  i+j\leq n} (\nu_i \cdot \nu_j')(-1),
\]
which by convention equals 1 for $n\leq 2$.

Define a meromorphic function 
\[
\gamma_\psi(s; \nu;\nu'):=\prod_{i+j\leq n}\gamma(s,  \nu_i\cdot \nu_j', \psi),
\]
and likewise an entire function 
\[
\varepsilon_\psi(s; \nu;\nu'):=\prod_{i+j\leq n}\varepsilon(s,  \nu_i\cdot \nu_j', \psi),
\]
which by convention is equal to 1 if $n=1$. Here $\gamma(s, \nu_i\cdot \nu_j',\psi)$ and $\varepsilon(s,  \nu_i\cdot \nu_j', \psi)$  are respectively the local  gamma factor and the local epsilon factor defined following the standard references \cite{T79, J79, K03}, which will be recalled below.

Given a  character $\omega$ of $\sk^\times$, %let $\oL(s,\omega)$ be the local $\oL$-factor of $\omega$. 
 Tate's local zeta integral (\cite{T50}) is defined by
\[
\oZ(s,\omega, \phi)=\int_{\sk^\times}\phi(x)\omega(x)\abs{x}_\sk^s\od^\times\!x,\quad \phi \in\CS(\sk),\]
which converges absolutely when $\Re(s)>-\mathrm{ex}(\omega)$. Here $\od^\times x:=\frac{\od x}{\abs{x}_\rk}$, which is a Haar measure of $\rk^\times$.  The local epsilon factor $\varepsilon(s,\omega,\psi)$ is an entire function defined by the local functional equation 
\begin{equation}\label{fe}
\frac{\oZ(1-s,\omega^{-1}, \widehat{\phi})}{\oL(1-s,\omega^{-1})}=\varepsilon(s,\omega,\psi)\cdot\frac{\oZ(s,\omega, \phi)}{\oL(s,\omega)},\quad  \phi \in \CS(\sk),
\end{equation}
where $\widehat{\phi}:=\CF_\psi(\phi) \in \CS(\sk)$ is the Fourier transform of $\phi$ with respect to $\psi$ defined by
\[
\CF_\psi(\phi)(x):=\int_{\sk}\phi(y)\psi(xy)\od\! y,\quad x\in \sk.
\]
% Recall that the Haar measure on $\K$ is self-dual with respect to $\psi$, so that it holds the Fourier inversion 
%\[
%\CF_{\psi^{-1}}\circ \CF_\psi(f)=f,\quad \forall f\in \CS(\sk).
%\] 
The meromorphic function 
\[
\gamma(s,\omega,\psi):=\varepsilon(s,\omega,\psi)\cdot\frac{\oL(1-s,\omega^{-1})}{\oL(s,\omega)}
\] 
is called the local gamma factor attached to $\omega$ and $\psi$.

\begin{remarkd} \label{rmk-gamma}
 A different convention is used in \cite{J09} by setting $\widehat{\phi}=\CF_{\overline{\psi}}(\phi)$ in \eqref{fe}. This changes the definition of $\varepsilon(s,\omega,\psi)$ by a factor $\omega(-1)$, thanks to the relation 
 \[
 \varepsilon(s,\omega,\overline{\psi})=\omega(-1)\cdot \varepsilon(s,\omega,\psi).
 \]
 We will translate the results in \cite{J09} according to the convention of this article. 
\end{remarkd}

Finally, define the meromorphic function 
\begin{equation} \label{Gamma}
\Gamma_\psi (s;\nu;\nu'):=\sgn(\nu;\nu')\cdot\gamma_\psi(s; \nu; \nu').
\end{equation} 
 Now we  state the main result of this paper.

%If  $n'=n$, put
%\[
 % \gamma_\psi(s; \nu;\nu'):=  \gamma_\psi(s; \nu; ({\nu'}_1,  {\nu'}_2, \cdots, {\nu'}_{n-1}))
%\]
%and
%\[
 % \varepsilon_\psi(s; \nu;\nu'):=  \varepsilon_\psi(s; \nu; ({\nu'}_1,  {\nu'}_2, \cdots, {\nu'}_{n-1}))
%\]

\begin{thmd} \label{thmA}
Assume that $s\in \Omega_{\nu, \nu'}$ as in Proposition \ref{prop1.2}. Let  $f\in I_\nu$, $f'\in I_{\nu'}$ and $\phi\in \CS(\sk^{1\times n})$.

\noindent 
(a)  If $n'=n$, then
\[
\Lambda(s, f, f', \phi)=\Gamma_\psi (s;\nu;\nu')\cdot \oZ(s, f, f',\phi).
\]

\noindent 
(b) If $n'=n-1$, then
\[
\Lambda(s, f, f')=\Gamma_\psi (s; \nu; \nu') \cdot \oZ(s, f, f').
\]

\end{thmd}

\begin{remarkd}
In Theorem \ref{thmA}, if $n'=n$ then
\[
\begin{aligned}
 & \ \Gamma_\psi (s;\nu;\nu')\cdot \oZ(s, f, f',\phi) \\
=& \  \sgn(\nu;\nu') \cdot \gamma_\psi(s; \nu; \nu') \cdot \oL(s, \nu\times \nu') \cdot \oZ^\circ(s, f, f',\phi)\\
=& \ \sgn(\nu;\nu') \cdot  \varepsilon_\psi(s;\nu;\nu')\cdot \prod_{i+j\leq n}\oL(1-s, \nu_i^{-1}\cdot \nu_j'^{-1})\cdot\prod_{i+j>n} \oL(s, \nu_i\cdot\nu_j') \cdot \oZ^\circ(s, f, f', \phi),
\end{aligned}
\]
which is  easily seen to be holomorphic in $s\in \Omega_{\nu, \nu'}$.

Likewise, if $n'=n-1$ then $\Gamma_\psi (s;\nu;\nu')\cdot \oZ(s, f, f')$ is holomorphic in $s\in \Omega_{\nu, \nu'}$ as well.
\end{remarkd}

When $\rk$ is Archimedean, $n'=n-1$, and $f$ and $f'$ lie in the minimal $K$-types (in the sense of Vogan), the Rankin-Selberg integrals $\oZ(s, f, f')$ have been explicitly calculated by Ishii and Miayzaki in \cite{IM22}.
They also obtain similar result for $n'=n$.  Still in the Archimedean case, the Rankin-Selberg integrals for minimal $K$-type vectors of irreducible generalized principal series representations  of $\GL(3)\times \GL(2)$ have been explicitly calculated by in Hirano, Ishii and Miyazaki in \cite{HIM22}. 

%Theorem \ref{thmA} in the Archimedean case is used in \cite{LLS21} to  prove the  Archimedean period relations for Rankin-Selberg convolutions for $\GL(n)\times \GL(n-1)$, which implies the period relations for critical values of the Rankin-Selberg L-functions for $\GL(n)\times \GL(n-1)$. 

This article is organized as follows. In Section \ref{sec2} we recall the Godement sections and  their basic properties. In Section \ref{sec3} we prove the recurrence relations for our integrals in terms of Godement sections. Proposition \ref{prop1.2} and Theorem \ref{thmA} will be proved in Sections \ref{sec4} and  \ref{sec5}  respectively, by using induction and the  recurrence relations. %,  Section is devoted to a proof of . The main result  by  for our integrals and the Rankin-Selberg integrals. 

\section{The Godement sections}\label{sec2}

We do not claim any originality of the results in this section. See \cite{J09} and \cite{IM22}. 
%However we will establish the analytic properties of Godement sections. 
Recall that $n$ is a positive integer, $n'=n$ or $n-1$, $\nu=(\nu_1, \nu_2, \cdots, \nu_n)\in (\widehat{\sk^\times})^n$ and $\nu'=(\nu'_1, \nu'_2, \cdots, \nu'_{n'})\in (\widehat{\sk^\times})^{n'}$.  Let $\chi\in \widehat{\sk^\times}$.

\subsection{A convergence result}

When $X$ is a Nash manifold or a totally disconnected locally compact topological space, we say that a function $f$ on $X$ is rapidly decreasing if 
\[
\abs{f(x)}\leq \phi(x),\quad \textrm{for all } x\in X,
\]
for some real valued function $\phi\in \CS(X)$.

\begin{lemd}\label{tate0}
Assume that 
\[
\mathrm{ex}(\nu_i')> i-1\quad\textrm{for all }  1\leq i\leq n'.
\]
Then the integral
\be\label{inttate}
   \int_{\bar B_{n'}} \phi(\bar b)\cdot \nu'(\bar b) \od\! \bar b
\ee
is absolutely convergent for all continuous functions  $\phi$ on $\sk^{n'\times n'}$ that are rapidly decreasing.
\end{lemd}
\begin{proof}
For every $a:=(a_{1},  a_{ 2}, \, \cdots,\,  a_{n'})\in \sk^{n'}$, put 
\[
 \phi_1(a)=\int_{\bar \n_{n'}} \abs{\phi(\bar u+\diag(a_{1},  a_{ 2}, \, \cdots,\,  a_{n'}))}\od\! \bar u,
\]
where $\bar \n_{n'}\subset \g\l_{n'}(\sk)$ is the subspace of the lower triangular nilpotent matrices, $\od\! \bar u$ is the product measure on $\bar \n_{n'}$ similar to \eqref{mea1}, and $\diag$ indicates the diagonal matrix. Then $\phi_1$ is a continuous function  on $\sk^{n'}$ that is rapidly decreasing in the archimedean case and has compact support in the non-archimedean case. 

Note that 
\begin{eqnarray*}
&& \int_{\bar B_{n'}} \abs{\phi(\bar b)\cdot \nu'(\bar b) }\od\! \bar b\\
&=&\int_{(\sk^\times)^{n'}}\phi_1(a)\cdot \abs{\nu'(a)} \cdot \prod_{i=1}^{n'} \absk{a_i}^{-i+1}\, \od^\times \! a,
\end{eqnarray*}
where $\od^\times \! a$ is the product of the Haar measures on $\sk^\times$. Then the lemma follows by 
 the usual argument in Tate's thesis. \end{proof}

\subsection{The Godement sections $G_{n-1}\longrightarrow G_n$}
In this subsection, we assume that $n'=n-1$.
For all $f'\in I_{\nu'}$ and $\phi\in \CS(\sk^{n'\times n})$, put
\be\label{gs10}
\begin{array}{rcl}
   \mathrm g^+(\nu',\chi, f', \phi):=  \int_{G_{n'}} f'(h^{-1})\cdot \phi([h,0])\cdot \chi(\det h)\cdot \absk{\det h}^{\frac{n}{2}}\od\! h.
\end{array}
\ee
%Put
%\[
 %\Omega_{n}^{(1)}:=\{(\nu'_1, \nu'_2, \cdots, \nu'_{n-1},\chi )\in (\widehat{\sk^\times})^{n}\,:\, \mathrm{ex}(\chi)>\mathrm{ex}(\nu_i')-1\quad\textrm{for all }  1\leq i\leq n-1\}.
%\]

\begin{prpd}\label{gsnnp}
Assume that 
\be\label{om1}
 \mathrm{ex}(\chi)>\mathrm{ex}(\nu_i')-1\quad\textrm{for all }  1\leq i\leq n-1.
\ee
Then the integral  \eqref{gs10} is  absolutely convergent.
\end{prpd}

\begin{proof}
This is proved in  \cite[Proposition 7.1 (i)]{J09}. 
\end{proof}

Suppose that \eqref{om1} holds. Then we have a well-defined map (see  \cite[Proposition 7.1 (iv)]{J09})
\be\label{gpp002}
\mathrm g^+_{\nu',\chi}:  I_{\nu'}\times \CS(\sk^{n'\times n}) \longrightarrow I_{(\nu',\chi)}
\ee
given by
\[
  (\mathrm g^+_{\nu',\chi}(f',\phi))(g):=\chi(\det g)\cdot \absk{\det g}^{\frac{n-1}{2}}\cdot \mathrm g^+(\nu',\chi,f', g. \phi),
\]
where $g\in G_n$ and $g.\phi$ indicates the right translation. 

\begin{prpd}\label{gssur1}
The image of the map \eqref{gpp002} spans the vector space  $I_{(\nu',\chi)}$.
%Suppose that $1$ is not a pole of the meromorphic function $\prod_{i=1}^{n'}\oL(s,{\nu_i'}^{-1}\cdot \chi)$. Then the image of the map \eqref{gpp} spans  $I_{(\nu',\chi)}$.
\end{prpd}
\begin{proof}
This directly follows from  \cite[Proposition 7.1 (v)]{J09}.
\end{proof}

\subsection{The Godement sections $G_{n}\longrightarrow G_n$}

For all $ f\in I_\nu$ and $\phi\in \CS(\sk^{n\times n})$,
put
\be\label{gs30}
\begin{array}{rcl}
   \mathrm g^\circ(\nu, \chi, f, \phi):= \int_{G_{n}} f(h) \cdot\phi(h)\cdot \chi(\det h)\cdot \abs{\det h}_{\rk}^{\frac{n-1}{2}}\od\! h.
\end{array}
\ee

\begin{prpd}\label{r2}
Assume that \be\label{conv1122}
  \mathrm{ex}(\chi)>-\mathrm{ex}(\nu_i)\quad\textrm{for all }\, i=1,2,\cdots, n.
\ee
Then the integral \eqref{gs30} is  absolutely convergent. 
\end{prpd}

\begin{proof}
The proof is similar to that of Proposition \ref{gsnnp}.
Fix a maximal compact subgroup $K_{n}$ of $G_{n}$, and fix a Haar measure on it such that
\be\label{dechg}
  \int_{G_{n}} \varphi(\bar b k)=\int_{\bar B_{n}} \int_{K_{n}}  \varphi(\bar b k)\od\! k\,\od \!\bar b
\ee
for all $\varphi\in \CS(G_{n})$.

Then we have that
\begin{eqnarray*}
&& \int_{G_{n'}}\left | f(h) \cdot\phi(h)\cdot \chi(\det h)\cdot \abs{\det h}^{\frac{n-1}{2}}_{\rk} \right |\od\! h\\
&=& \int_{\bar B_{n}} \phi_1(\bar b)\cdot \abs{\nu(\bar b)}\cdot  \abs{\chi(\det \bar b)}\cdot \absk{\det \bar b}^{\frac{n-1}{2}}\cdot  \bar \rho_{n}(\bar b)\od \! \bar b, 
\end{eqnarray*}
where 
\[
 \phi_1(\bar b)=\int_{K_{n}} \left| f( k)\cdot \phi(\bar b k)\right|\od\! k.
\]
The Proposition  then follows from Lemma \ref{tate0}.

\end{proof}

Assume that  \eqref{conv1122} holds. Define a map
\be\label{gpp2}
\mathrm g^\circ_{\nu,\chi}:  I_{\nu}\times \CS(\sk^{n\times n}) \longrightarrow I_{\nu}
\ee
by
\[
  (\mathrm g^\circ_{\nu,\chi}(f,\phi))(g):=\chi(\det g^{-1})\cdot \absk{\det g}^{\frac{1-n}{2}}\cdot \mathrm g^\circ(\nu,\chi,f, L_g \phi),
\]
where $g\in G_n$, and $L_g$  stands for the left translation so that $(L_g\phi)(x)=\phi(g^{-1} x)$ for all $x\in \sk^{n\times n}$. It is easy to see that this map is well-defined and bilinear (\cf \cite[Proposition 3.2]{IM22}).

\begin{prpd}
Assume the condition \eqref{conv1122}.   Then the image of the map \eqref{gpp2} spans  $I_{\nu}$.
\end{prpd}
\begin{proof}
By change of variable, we have that
\[
  (\mathrm g^\circ_{\nu,\chi}(f, \phi))(g)
 =\int_{G_{n}} \phi(h)\cdot f(gh)\cdot \chi(\det h)\cdot \abs{\det h}^{\frac{n-1}{2}}_{\rk} \od\! h.
\]
Note that $\CS(G_n)\subset \CS(\sk^{n\times n})$, and then the proposition easily follows by Dixmier-Malliavin Theorem \cite{DM78}.
\end{proof}

\section{Recurrence relations} \label{sec3}

We continue with the notation of the last section. Recall that $\chi\in \widehat{\sk^\times}$.
For $s\in \C$,  write $\chi_s:=\chi\cdot \abs{\,\cdot\,}_\sk^s\in \widehat{ \sk^\times}$. 
Suppose that $n'=n-1$ in this section, so that $\nu\in (\widehat{\rk^\times})^n$ and $\nu' \in (\widehat{\rk^\times})^{n-1}$.

\subsection{The first recurrence relation}

\begin{prpd}\label{r1}
Let $\phi_1\in \CS(\sk^{(n-1)\times n})$ and $\phi_2\in \CS(\sk^{1\times n})$, and write $\phi_0:=\phi_1\otimes \phi_2\in  \CS(\sk^{n\times n})$. Then for all $f\in I_\nu$ and $f'\in I_{\nu'}$, 
\[
   \Lambda(s, f,    \mathrm g_{\nu',\chi}^{+}(f',\phi_1), \phi_2)=\Lambda(s,\mathrm g^\circ_{\nu,\chi_s}(f, \phi_0), f').
\]
 %More precisely, the left hand side double integral is absolutely convergent if and only if so is the right one, and when this is the case the two integrals are equal to each other.
\end{prpd}

As explained right below Lemma \ref{lem1.1}, the equation in Proposition \ref{r1} should be understood that both sides have the same range of absolute convergence. This will be clear from the proof below.

\begin{proof}
We have that
\begin{eqnarray}
\label{lc0}
&& \Lambda(s,\mathrm g^\circ_{\nu,\chi_s}(f, \phi_0)), f')\\
\nonumber &=& \int_{G_{n-1}} \mathrm g^\circ_{\nu,\chi_s}(f, \phi_0)\left(z_n \left[
                                                                \begin{array}{cc}
                                                                  h & 0 \\
                                                                  0 & 1 \\
                                                                \end{array}
                                                              \right]
 \right) f'(z_{n-1} h) \cdot \abs{\det h}^{s-\frac{1}{2}}_{\rk} \od\! h\\
\nonumber &=& \int_{G_{n-1}} \int_{G_{n}}  \phi_0(g)\cdot f\left(z_n \left[
                                                                \begin{array}{cc}
                                                                  h & 0 \\
                                                                  0 & 1 \\
                                                                \end{array}
                                                              \right] g\right)\cdot \chi(\det g)\cdot \abs{\det g}^{s+\frac{n-1}{2}}_{\rk} \od\! g \\
               \nonumber                                               &&\qquad \cdot f'(z_{n-1} h) \cdot \abs{\det h}^{s-\frac{1}{2}}_{\rk} \od\! h.
                                                           \end{eqnarray}
                By the change of variables $g\longmapsto  \left[
                                                                \begin{array}{cc}
                                                                  h^{-1} & 0 \\
                                                                  0 & 1 \\
                                                                \end{array}
                                                              \right]g $, the above inner integral equals
        \[
         \chi(\det h^{-1})\cdot \abs{\det h}^{-s-\frac{n-1}{2}}_{\rk} \cdot \int_{G_{n}}  \phi_0\left(\left[
                                                                \begin{array}{cc}
                                                                  h^{-1} & 0 \\
                                                                  0 & 1 \\
                                                                \end{array}
                                                              \right]g\right)\cdot f\left(z_n  g\right)\cdot \chi(\det g)\cdot \abs{\det g}^{s+\frac{n-1}{2}}_{\rk} \od\! g. \\
        \]
        Hence
        \begin{eqnarray*}
        \eqref{lc0} &=&     \int_{G_{n-1}}   \int_{G_{n}}  \phi_0\left(\left[
                                                                \begin{array}{cc}
                                                                  h^{-1} & 0 \\
                                                                  0 & 1 \\
                                                                \end{array}
                                                              \right]g\right)\cdot f\left(z_n  g\right)\cdot \chi(\det g)\cdot \abs{\det g}^{s+\frac{n-1}{2}}_{\rk} \od\! g \\
                                                              &&\qquad \cdot f'(z_{n-1} h) \cdot \chi(\det h^{-1})\cdot \abs{\det h}^{-\frac{n}{2}}_{\rk} \od\! h.\\
         \end{eqnarray*}
         We are free to change the order of the integrals since we are only concerned with the case when the double integrals are absolutely convergent. Thus we have that
    \begin{eqnarray*}
        \eqref{lc0} &=&    \int_{G_{n}}   \int_{G_{n-1}}    \phi_0\left(\left[
                                                                \begin{array}{cc}
                                                                  h^{-1} & 0 \\
                                                                  0 & 1 \\
                                                                \end{array}
                                                              \right]g\right) \cdot f'(z_{n-1} h) \cdot \chi(\det h^{-1})\cdot \abs{\det h}^{-\frac{n}{2}}_{\rk} \od\! h\\
                                                              &&\qquad  \cdot f\left(z_n  g\right)  \cdot \chi(\det g)\cdot \abs{\det g}^{s+\frac{n-1}{2}}_{\rk} \od\! g \\
                                                              &=&    \int_{G_{n}}   \int_{G_{n-1}}    \phi_1\left(\left[
                                                                \begin{array}{cc}
                                                                  h^{-1} & 0 \\
                                                               \end{array}
                                                              \right]g\right) \cdot f'(z_{n-1} h) \cdot \chi(\det h^{-1})\cdot \abs{\det h}^{-\frac{n}{2}}_{\rk} \od\! h\\
                                                              &&\qquad  \cdot f\left(z_n  g\right) \cdot\phi_2(e_n g) \cdot \chi(\det g)\cdot \abs{\det g}^{s+\frac{n-1}{2}}_{\rk} \od\! g. \\
         \end{eqnarray*}
         By the change of variables $h\longmapsto z_{n-1}^{-1} h^{-1}$, the above inner integral equals
    \begin{eqnarray*}
          &&  \int_{G_{n-1}}    \phi_1\left(\left[
                                                                \begin{array}{cc}
                                                                  h & 0 \\
                                                               \end{array}
                                                              \right]  \left[
                                                                \begin{array}{cc}
                                                                  z_{n-1} & 0 \\
                                                                  0 & 1 \\
                                                                \end{array}
                                                              \right] g\right) \cdot f'(h^{-1}) \cdot \chi(\det h)\cdot \abs{\det h}^{\frac{n}{2}}_{\rk} \od\! h\\
                                                              &=&\chi(\det g^{-1})\cdot \abs{\det g}^{-\frac{n-1}{2}}_{\rk} \cdot (\mathrm g^+_{\nu',\chi}(f', \phi_1))\left( \left[
                                                                \begin{array}{cc}
                                                                  z_{n-1} & 0 \\
                                                                  0 & 1 \\
                                                                \end{array}
                                                              \right] g\right).
      \end{eqnarray*}
         Therefore
          \begin{eqnarray*}
        \eqref{lc0} &=&    \int_{G_{n}}     f\left(z_n  g\right)  \cdot (\mathrm g^+_{\nu',\chi}(f', \phi_1))\left( \left[
                                                                \begin{array}{cc}
                                                                  z_{n-1} & 0 \\
                                                                  0 & 1 \\
                                                                \end{array}
                                                              \right] g\right) \cdot  \phi_2(e_n g)\cdot   \abs{\det g}^{s}_{\rk} \od\! g \\
                                                              &=&  \Lambda(s, f,    \mathrm g_{\nu',\chi}^{+}(f',\phi_1), \phi_2).
         \end{eqnarray*}
         This finishes the proof of the proposition.
\end{proof}

\subsection{The second recurrence relation} 
For every $\phi\in \CS(\sk^{k\times l})$ ($k,l\in \BN$), write ${}^t\phi\in \CS(\sk^{l\times k})$ for the function 
\[
{}^t\phi(x)= \phi({}^tx).
\]
 For every $f\in \con^\infty(G_k)$, write $\widehat {f}\in  \con^\infty(G_k)$ for the function 
 \[
\widehat {f}(g) = f(w_k g^\iota w_k )
 \]
 where $g^\iota:={}^t g^{-1}$ is the transpose inverse of $g$. For every  $\alpha=(\alpha_1, \alpha_2, \cdots, \alpha_k)\in (\widehat{\sk^\times})^k$, write
\[
  \widehat \alpha:=(\alpha_k^{-1}, \alpha_{k-1}^{-1}, \cdots, \alpha_1^{-1})\in (\widehat{\sk^\times})^k.
\]
Then $\widehat f\in I_{\widehat \alpha}$ whenever $f\in I_\alpha$.

\begin{prpd}\label{prprl1}
Let
  $\mu\in (\widehat{\sk^\times})^{n-1}$.
Let $\phi_1\in \CS(\sk^{(n-1)\times (n-1)})$ and $\phi_2\in \CS(\sk^{(n-1)\times 1})$, and write $\phi_0:=\phi_1\otimes \phi_2\in  \CS(\sk^{(n-1)\times n})$. Then for all $f_\mu \in I_\mu$ and $f' \in I_{\nu'}$, 
\[
  \Lambda(s,\mathrm g^+_{\mu, \chi}(f_\mu, \phi_0), f')=\Lambda(1-s,\widehat{\mathrm g^\circ_{\nu', \chi_s}(f' , \phi_1)}, \widehat{f_\mu}, w_{n-1}.{}^t\phi_2).
\]
% More precisely, the left hand side double integral is absolutely convergent if and only if so is the right one, and when this is the case the two integrals equal to each other.
\end{prpd}

\begin{proof}
The proposition is trivial when $n=1$. Thus we assume that $n\geq 2$.
We have that
\begin{eqnarray}
\label{lc02}&& \Lambda(1-s,\widehat{\mathrm g^\circ_{\nu', \chi_s}(f', \phi_1)}, \widehat{ f_{\mu}}, w_{n-1}.{}^t\phi_2)\\
\nonumber
 &=& \int_{G_{n-1}} \abs{\det g}^{1-s}_{\rk} \cdot  (\mathrm g^\circ_{\nu', \chi_s}(f', \phi_1))(w_{n-1} z_{n-1}^\iota g^\iota w_{n-1}) \\
 \nonumber  &&\qquad
                                                     \cdot f_{\mu}\left( w_{n-1}\left[
                                                                \begin{array}{cc}
                                                                  z_{n-2}^\iota & 0 \\
                                                                  0 & 1 \\
                                                                \end{array}
                                                              \right] g^\iota w_{n-1}\right) \cdot {}^t\phi_2(e_{n-1} g w_{n-1})
                                                              \od\! g \\
                      \nonumber    &=& \int_{G_{n-1}} \abs{\det g}^{1-s}_{\rk} \cdot \phi_2(w_{n-1} {}^tg {}^te_{n-1} )
                                                      \cdot f_{\mu}\left( w_{n-1}\left[
                                                                \begin{array}{cc}
                                                                  z_{n-2}^{\iota} & 0 \\
                                                                  0 & 1 \\
                                                                \end{array}
                                                              \right] g^{\iota} w_{n-1}\right)\\
                              \nonumber                              &&\qquad \cdot \int_{G_{n-1}}  \phi_1(h)\cdot f'(w_{n-1} z_{n-1}^{\iota} g^{\iota} w_{n-1} h)\cdot \chi(\det h)\cdot \abs{\det h}^{s+\frac{n-2}{2}}_{\rk} \od\! h
                                                              \od\! g.
 \end{eqnarray}
 By the change of variables $h\longmapsto w_{n-1} {}^tg {}^tz_{n-1} w_{n-1} z_{n-1} h$, the above inner integral equals
 \[
     \chi(\det g)\cdot \abs{\det g}^{s+\frac{n-2}{2}}_{\rk} \cdot \int_{G_{n-1}}  \phi_1( w_{n-1} {}^tg {}^tz_{n-1} w_{n-1} z_{n-1} h)\cdot f'( z_{n-1} h)\cdot \chi(\det h)\cdot \abs{\det h}^{s+\frac{n-2}{2}}_{\rk} \od\! h.
 \]
 As in the proof of Proposition \ref{r1},  we are free to change the order of the integrals.
Thus we have that
\begin{eqnarray*}
\eqref{lc02}&=&\int_{G_{n-1}}\chi(\det g)\cdot  \abs{\det g}^{\frac{n}{2}}_{\rk} \cdot \phi_2(w_{n-1} {}^tg {}^te_{n-1} )
                                                      \cdot f_{\mu}\left( w_{n-1}\left[
                                                                \begin{array}{cc}
                                                                  z_{n-2}^{\iota} & 0 \\
                                                                  0 & 1 \\
                                                                \end{array}
                                                              \right] g^{\iota} w_{n-1}\right)\\
                                                 &&\qquad \cdot \int_{G_{n-1}}  \phi_1( w_{n-1} {}^tg {}^t z_{n-1} w_{n-1} z_{n-1} h)\cdot f'( z_{n-1} h)\cdot \chi(\det h)\cdot \abs{\det h}^{s+\frac{n-2}{2}}_{\rk} \od\! h      \od\! g\\
                                                                     &=&\int_{G_{n-1}}  f'( z_{n-1} h)\cdot \chi(\det h)\cdot \abs{\det h}^{s+\frac{n-2}{2}}_{\rk} \cdot \xi(h)
                                                    \od\! h,
\end{eqnarray*}
where
\begin{eqnarray*}
\xi(h)&:=&\int_{G_{n-1}}\chi(\det g)\cdot  \abs{\det g}^{\frac{n}{2}}_{\rk} \cdot  \phi_1( w_{n-1} {}^tg {}^tz_{n-1} w_{n-1} z_{n-1} h)
                                                      \cdot \phi_2(w_{n-1} {}^tg {}^te_{n-1} ) \\
                                                 &&\qquad \cdot f_{\mu}\left( w_{n-1}\left[
                                                                \begin{array}{cc}
                                                                  z_{n-2}^{\iota} & 0 \\
                                                                  0 & 1 \\
                                                                \end{array}
                                                              \right] g^{\iota} w_{n-1}\right)  \od\! g.
\end{eqnarray*}
By the change of variable $g\longmapsto   \left[
                                                                \begin{array}{cc}
                                                                  z_{n-2}^{-1} & 0 \\
                                                                  0 & 1 \\
                                                                \end{array}
                                                              \right] w_{n-1} {}^tg w_{n-1} $, we have that
\begin{eqnarray*}
\xi(h)& =&\int_{G_{n-1}}\chi(\det g)\cdot  \abs{\det g}^{\frac{n}{2}}_{\rk} \cdot  \phi_1\left(  g w_{n-1}\left[ \begin{array}{cc}
                                                                  z_{n-2}^{\iota} & 0 \\
                                                                  0 & 1 \\
                                                                \end{array}
                                                              \right] {}^tz_{n-1} w_{n-1} z_{n-1} h\right)
                                                       \\
                                                 &&\qquad \cdot \phi_2\left( g w_{n-1}\left[ \begin{array}{cc}
                                                                  z_{n-2}^{\iota} & 0 \\
                                                                  0 & 1 \\
                                                                \end{array}
                                                              \right]  {}^te_{n-1} \right) \cdot f_{\mu}( g^{-1} )  \od\! g.\\
                                                            \end{eqnarray*}
Recall that
\[
       z_n :=
 \left[
            \begin{array}{cc}
                            w_{n-1}& 0 \\
                       0 & 1 \\
                     \end{array}
                   \right]
                    \left[
            \begin{array}{cc}
                       z_{n-2}^{\iota}& 0 \\
                       0 & 1_2 \\
                     \end{array}
                   \right]
                    \left[
            \begin{array}{cc}
                       {}^tz_{n-1} w_{n-1} z_{n-1}& {}^te_{n-1}\\
                       0 & 1 \\
                     \end{array}
                   \right].
                   \]
This implies that
\[
  [g,0] z_n \left[ \begin{array}{cc}
                                                                  h & 0 \\
                                                                  0 & 1 \\
                                                                \end{array}
                                                              \right] =\left[g w_{n-1}\left[ \begin{array}{cc}
                                                                  z_{n-2}^{\iota} & 0 \\
                                                                  0 & 1 \\
                                                                \end{array}
                                                              \right] {}^t z_{n-1} w_{n-1} z_{n-1} h , g w_{n-1} \left[ \begin{array}{cc}
                                                                  z_{n-2}^{\iota} & 0 \\
                                                                  0 & 1 \\
                                                                \end{array}
                                                              \right] {}^t e_{n-1} \right],
\]
 which further implies that
 \begin{eqnarray*}
 && \phi_1\left(  g w_{n-1}\left[ \begin{array}{cc}
                                                                  z_{n-2}^{\iota} & 0 \\
                                                                  0 & 1 \\
                                                                \end{array}
                                                              \right] {}^t z_{n-1} w_{n-1} z_{n-1} h\right)
                                                         \cdot \phi_2\left( g w_{n-1}\left[ \begin{array}{cc}
                                                                  z_{n-2}^{\iota} & 0 \\
                                                                  0 & 1 \\
                                                                  \end{array}
                                                                  \right] {}^te_{n-1}\right)\\
                                                                  &=& \phi_0\left( [g,0] z_n \left[ \begin{array}{cc}
                                                                  h & 0 \\
                                                                  0 & 1 \\
                                                                \end{array}
                                                              \right]\right).
                                                                  \end{eqnarray*}
 Hence
 \[
 \xi(h)=\int_{G_{n-1}}\chi(\det g)\cdot  \abs{\det g}^{\frac{n}{2}}_{\rk} \cdot \phi_0\left( [g,0] z_n \left[ \begin{array}{cc}
                                                                  h & 0 \\
                                                                  0 & 1 \\
                                                                \end{array}
                                                              \right]\right) \cdot f_{\mu}( g^{-1} )  \od\! g,\\
                                                              \]
  and
  \[
  \chi(\det h)\cdot \abs{\det h}^{\frac{n-1}{2}}_{\rk} \cdot \xi(h)= (\mathrm g^+_{\mu, \chi}(f_{\mu}, \phi_0))\left(z_n \left[ \begin{array}{cc}
                                                                  h & 0 \\
                                                                  0 & 1 \\
                                                                \end{array}\right]\right).
  \]

Finally,
\begin{eqnarray*}
\eqref{lc02}&=&\int_{G_{n-1}}  f'( z_{n-1} h)\cdot \abs{\det h}^{s-\frac{1}{2}}_{\rk} \cdot (\mathrm g^+_{\mu,\chi}(f_{\mu}, \phi_0))\left(z_n \left[ \begin{array}{cc}
                                                                  h & 0 \\
                                                                  0 & 1 \\
                                                                \end{array}\right]\right)
                                                    \od\! h\\
                                                 &=&     \Lambda(s,\mathrm g^+_{\mu,\chi}(f_{\mu}, \phi_0)), f').
                                                    \end{eqnarray*}
  This finishes the proof of the proposition.
\end{proof}

\section{Absolute convergence} \label{sec4}

In this section we prove Proposition \ref{prop1.2}, which is restated below. Let $f\in I_\nu$, $f'\in I_{\nu'}$ and $\phi\in \CS(\sk^{1\times n})$, as in Theorem \ref{thmA}.

\begin{prpd}\label{prpconv1}

Assume that $s$ lies in 
\[
\Omega_{\nu, \nu'}:=\left\{ s\in \mathbb{C}\left|\begin{aligned} & \mathrm{ex}(\nu_i)+\mathrm{ex}(\nu'_j)+\Re(s)<1  \textrm{ whenever } i+j\leq n,\\
& \mathrm{ex}(\nu_i)+\mathrm{ex}(\nu'_j)+\Re(s)>0 \textrm{ whenever } i+j> n
\end{aligned}\right.\right\}.
\]
\noindent {\rm (a)} Assume that $n'=n$. Let $\mu':=(\nu'_1, \nu'_2, \cdots, \nu'_{n-1})$, $f_{\mu'}\in I_{\mu'}$, and $\phi_1\in \CS(\sk^{(n-1)\times n})$. Then the double integral 
\[
\Lambda(s,\mathrm g^\circ_{\nu,\nu'_n\cdot\absk{\,\cdot\,}^s}(f, \phi_0), f_{\mu'}),
\]
where $\phi_0:=\phi_1\otimes \phi\in \CS(\sk^{n\times n})$, is absolutely convergent. 

\noindent {\rm (b)} Assume that $n'=n-1$. Let  $\mu:=(\nu_1, \nu_2, \cdots, \nu_{n-1})$, $f_{\mu}\in I_{\mu}$, $\phi_2\in \CS(\sk^{(n-1)\times (n-1)})$, and $\phi_3\in \CS(\sk^{(n-1)\times 1})$. Then the double integral
\[
\Lambda(1-s,\widehat{\mathrm g^\circ_{\nu', \nu_n\cdot\absk{\,\cdot\,}^s}(f', \phi_2)}, \widehat{ f_{\mu}}, w_{n-1}.{}^t\phi_3)
\] 
is absolutely convergent. 

\noindent {\rm (c)} 
The integral $\Lambda(s, f, f', \phi)$ is absolutely convergent if $n'=n$, and the integral $\Lambda(s, f, f')$ is absolutely convergent if $n'=n-1$. 

\end{prpd}

\begin{proof}

We prove (c) by induction on $n+n'$. The statement (c) is trivial when $n+n'=1$. So we assume that $n+n'\geq 2$, and that (c) holds when $n+n'$ is smaller.
Without loss of generality, we assume that all the $\nu_i$'s and $\nu'_j$'s are positive characters.

We first assume that $n'=n$. The assumptions of the proposition imply that
\[
  \mathrm{ex}(\nu'_n)>\mathrm{ex}(\nu_i')-1, \quad \textrm{for all } i=1,2, \cdots, n-1.
\]
In view of Proposition \ref{gssur1}, we assume without loss of generality that
\[
  f'=\mathrm g_{\mu',\nu'_n}^{+}(f_{\mu'},\phi_1),
\]
where $\mu'$, $f_{\mu'}$ and $\phi_1$ are as in (a).
By Proposition \ref{r1}, we have that
\be\label{rnn1}
   \Lambda(s, f,    f', \phi)=\Lambda(s,\mathrm g^\circ_{\nu,\nu'_n\cdot\absk{\,\cdot\,}^s}(f, \phi_0), f_{\mu'}),
\ee
where $\phi_0:=\phi_1\otimes \phi\in \CS(\sk^{n\times n})$.

Recall that \[
  (g^\circ_{\nu,\nu'_n\cdot\absk{\,\cdot\,}^s}(f, \phi_0))(g)
 =\int_{G_{n}} \phi_0(h)\cdot f(gh)\cdot \nu'_n(\det h)\cdot \absk{\det h}^{s+\frac{n-1}{2}}\od\! h,
\]
where $g\in G_n$.
This is absolutely convergent by Proposition \ref{r2}. Note that
the function
\[
  g\longmapsto \int_{G_{n}} \left|\phi_0(h)\cdot f(gh)\cdot \nu'_n(\det h)\cdot \absk{\det h}^{s+\frac{n-1}{2}}\right|\od \! h
\]
is bounded by a positive function in $I_{\nu}$. Thus the integral \eqref{rnn1} is absolutely convergent by the induction hypothesis.

Now we assume that $n'=n-1$. The assumptions of the proposition imply that
\[
  \mathrm{ex}(\nu_n)>\mathrm{ex}(\nu_i)-1, \quad \textrm{for all } i=1,2, \cdots, n-1.
\]
Proposition \ref{gssur1} implies that $\abs{f}$ is bounded by a finite sum of functions of the form
\[
  \abs{\mathrm g_{\mu,\nu_n}^{+}(f_{\mu},\phi_2\otimes \phi_3)},
\]
where $\mu$, $f_{\mu}$, $\phi_2$ and $\phi_3$ are as in (b), so that $\phi_2\otimes \phi_3\in \CS(\sk^{(n-1)\times n})$. Thus we may assume without loss of generality that
\[
 f=\mathrm g_{\mu,\nu_n}^{+}(f_{\mu},\phi_2\otimes \phi_3).
\]
Then by Proposition \ref{prprl1}, we have that
\be\label{rnn12}
  \Lambda(s,f, f')=\Lambda(1-s,\widehat{\mathrm g^\circ_{\nu', \nu_n\cdot\absk{\,\cdot\,}^s}(f', \phi_2)}, \widehat{ f_{\mu}}, w_{n-1}.{}^t\phi_3).
\ee

Recall that \[
  (g^\circ_{\nu',\nu_n\cdot\absk{\,\cdot\,}^s}(f', \phi_2))(g)
 =\int_{G_{n-1}} \phi_2(h)\cdot f'(gh)\cdot \nu_n(\det h)\cdot \absk{\det h}^{s-1+\frac{n}{2}}\od \! h,
\]
where $g\in G_{n-1}$.
This is absolutely convergent by Proposition \ref{r2}. Note that
the function
\[
  g\longmapsto   \int_{G_{n-1}} \left| \phi_2(h)\cdot f'(w_{n-1}g^{\iota} w_{n-1} h)\cdot \nu_n(\det h)\cdot \absk{\det h}^{s-1+\frac{n}{2}}\right| \od\! h
\]
is bounded by a positive function in $I_{\widehat{\nu'}}$. It is elementary to check that
\[
  \mathrm{ex}({\nu'}_{n-i}^{-1})+\mathrm{ex}(\nu_{n-j}^{-1})+1-\Re(s)<1 \quad \textrm{for all } 1\leq i,j\leq n-1 \textrm{ and } i+j\leq n-1,
\]
and
\[
 \mathrm{ex}({\nu'}_{n-i}^{-1})+\mathrm{ex}(\nu_{n-j}^{-1})+1-\Re(s)>0 \quad \textrm{for all } 1\leq i,j\leq n-1 \textrm{ and } i+j> n-1.
\]
Thus the integral \eqref{rnn12} is absolutely convergent by the induction hypothesis.

This finishes the proof of (c). As we mentioned, the function 
\[
  g\longmapsto \int_{G_{n}} \left|\phi_0(h)\cdot f(gh)\cdot \nu'_n(\det h)\cdot \absk{\det h}^{s+\frac{n-1}{2}}\right|\od \! h
\]
is bounded by a positive function in $I_{\nu}$. Together with (c), this proves (a). Similarly (b) also follows from (c).
\end{proof}

\section{Proof of Theorem \ref{thmA}} \label{sec5}

In this section we prove Theorem \ref{thmA} by induction on $n+n'$. As before, let $f\in I_\nu$, $f'\in I_{\nu'}$ and $\phi\in \CS(\sk^{1\times n})$.

\begin{prp}\label{prop5.1}
If $n=1$, then 
\[
\Lambda(s, f, f', \phi)=\oZ(s, f, f', \phi).
\]
\end{prp}

\begin{proof}
This is straightforward from the definitions of the two integrals. 
\end{proof}

Proposition \ref{prop5.1} implies that Theorem \ref{thmA} holds when $n+n'\leq 2$. In the rest of this subsection we assume that $n+n'\geq 3$.

\begin{prp} \label{prop5.2} Assume that Theorem \ref{thmA} holds for $G_n\times G_{n-1}$. Then it holds for $G_n\times G_n$.

\end{prp}

\begin{proof}
Suppose that $n'=n$. Write $\mu'=(\nu_1',\ldots, \nu_{n-1}')$ so that $\nu'=(\mu', \nu_n')$. Note that
\[
 \Gamma_\psi (s; \nu;\nu')= \Gamma_\psi (s; \nu;\mu').
\]

By Proposition \ref{r1} we have that
\be\label{doublei2}
 \Lambda(s, f,     \mathrm g^{+}_{\mu', \nu_n'}(f'_{\mu'},\phi_1), \phi) \\
= \Lambda(s,\mathrm g^\circ_{\nu, \nu_n' \cdot\abs{\cdot}_{\sk}^s}(f, \phi_0)), f'_{\mu'}),
\ee
where $f'_{\mu'}\in I_{\mu'}$, $\phi_1\in \CS(\sk^{(n-1)\times n})$, and $\phi_0:=\phi_1\otimes \phi\in  \CS(\sk^{n\times n})$. Since $s\in \Omega_{\nu, \nu'}$,  Proposition \ref{prpconv1} (a) implies that the double integrals in both sides of \eqref{doublei2} are absolutely convergent. 

On the other hand, by \cite[Proposition 3.4]{IM22},
\[
 \oZ(s, f,     \mathrm g^{+}_{\mu', \nu_n'}(f'_{\mu'},\phi_1), \phi) \\
= \oZ(s,\mathrm g^\circ_{\nu, \nu_n' \cdot\abs{\cdot}_{\sk}^s}(f, \phi_0)), f'_{\mu'}).
\]
In view of Proposition \ref{gssur1}, the proposition follows from the above two equalities. 
\end{proof}

\begin{prp} \label{prop5.3} Assume that Theorem \ref{thmA} holds for $G_{n-1}\times G_{n-1}$. Then it holds for $G_n\times G_{n-1}$.
%and \begin{eqnarray*}
%&&  \Gamma(s;\nu; \nu')=\prod_{1\leq i\leq n-1}  \nu'_i(-1)^n \cdot \prod_{1\leq i, j \leq n-1} \gamma(s, \nu_i\cdot \nu'_j, \psi) \\
% &&\qquad \cdot  \Gamma(1-s; ({\nu'}_{n-1}^{-1},  {\nu'}_{n-2}^{-1}, \cdots, {\nu'}_{1}^{-1}) ;  ({\nu}_{n-1}^{-1},  {\nu}_{n-2}^{-1}, \cdots, {\nu}_{1}^{-1})).
% \end{eqnarray*}
\end{prp}

\begin{proof}
Suppose that $n'=n-1$. Write $\mu=(\nu_1,\ldots, \nu_{n-1})$ so that $\nu=(\mu, \nu_n)$. By \cite[Proposition 3.5]{IM22}, 
\begin{equation} \label{5.3eq1}
\oZ(s, \mathrm{g}^+_{\mu, \nu_n}(f_\mu, \phi_0), f')=\oZ(s, f_\mu, \mathrm{g}^\circ_{\nu', \nu_n\cdot\abs{\cdot}^s_\sk}(f', \phi_1), \CF_{\overline{\psi}}(\phi_2)),
\end{equation}
where $f_\mu\in I_\mu$ and $\phi_0=\phi_1\otimes\phi_2\in \CS(\sk^{(n-1)\times n})$ are as in Proposition \ref{prprl1}, and $\CF_{\overline{\psi}}(\phi_2)\in \CS(\sk^{1\times(n-1)})$ is the Fourier transform of $
\phi_2\in \CS(\sk^{(n-1)\times 1})$ with respect to $\overline{\psi}$ defined by
\[
\CF_{\overline{\psi}}(\phi_2)(x):=\int_{\sk^{(n-1)\times 1}}\phi_2(y)\overline{\psi}(xy) \od\! y,\quad x \in \sk^{1\times (n-1)}.
\] 

We now apply the functional equation of Rankin-Selberg integrals as in \cite{J09}, where the notations are slightly different from ours. %For $W\in C^\infty(G_k)$ write $\widetilde{W}\in C^\infty(G_k)$ for the function 
%$
%\widetilde{W}(g)= W(w_k g^\iota).
%$
%It is straightforward to check that 
%\[
%\widetilde{W_f}(g)=w_k. \overline{W}_{\widehat{f}}(g).
%\]
Put
\[
\widetilde{\phi_2}(x):=(w_{n-1}.{}^t\phi_2)(-x), \quad x\in \sk^{1\times(n-1)}.
\]
Then $\widetilde{\phi_2}\in \CS(\sk^{1\times (n-1)})$. Write
\[
   \gamma(s, I_\mu\times I_{\nu'}, \psi):=\prod_{1\leq i,j\leq n-1}\gamma(s,  \nu_i\cdot \nu_j', \psi),
\]
which is a product of local gamma factors. 

By \cite[Theorem 2.1]{J09} and Remark \ref{rmk-gamma}, and by noting that 
\[  ((\CF_{\overline{\psi}}\circ\CF_{\overline{\psi}})(\phi_2))(x)=\phi_2(-x), \quad x\in \sk^{(n-1)\times 1},  \]
 we obtain that
\begin{equation} \label{5.3eq2}
\begin{aligned}
& \omega_{\mu}(-1)\cdot \omega_{\nu'}(-1)^{n-1} \cdot  \gamma(s, I_\mu\times I_{\nu'}, \psi)\cdot \oZ(s, f_\mu, \mathrm{g}^\circ_{\nu', \nu_n\cdot\abs{\cdot}^s_\sk}(f', \phi_1), \CF_{\overline{\psi}}(\phi_2)) \\ 
= \ & \oZ(1-s,  \widehat{\mathrm g^\circ_{\nu', \nu_n\cdot\abs{\cdot}^s_\sk}(f', \phi_1)} , \widehat f_\mu, \widetilde{\phi_2}),
 \end{aligned}
\end{equation}
where $\omega_{\mu}$ and $\omega_{\nu'}$ denote the central characters of $I_{\mu}$ and $I_{\nu'}$ respectively, so that
\[
\omega_{\mu}(-1)=\prod_{1\leq i\leq n-1}\nu_i(-1),\quad \omega_{\nu'}(-1)=\prod_{1\leq i\leq n-1}\nu'_i(-1).
\]

Combining \eqref{5.3eq1} and \eqref{5.3eq2}, we obtain that
\be \label{5.3eq1+2}
\begin{aligned}
& \omega_{\mu}(-1)\cdot \omega_{\nu'}(-1)^{n-1} \cdot  \gamma(s, I_\mu\times I_{\nu'}, \psi)\cdot \oZ(s, \mathrm{g}^+_{\mu, \nu_n}(f_\mu, \phi_0), f') \\ 
= \ & \oZ(1-s,  \widehat{\mathrm g^\circ_{\nu', \nu_n\cdot\abs{\cdot}^s_\sk}(f', \phi_1)} , \widehat f_\mu, \widetilde{\phi_2}).
 \end{aligned}
\ee

By Proposition \ref{prprl1}, we have that
\begin{equation} \label{5.3eq3}
\begin{aligned}
&  \ \Lambda(s, \mathrm{g}^+_{\mu, \nu_n}(f_\mu, \phi_0), f')\\
 =&\ \Lambda(1-s,  \widehat{\mathrm g^\circ_{\nu', \nu_n\cdot\abs{\cdot}^s_\sk}(f', \phi_1)} , \widehat f_\mu, w_{n-1}.{}^t\phi_2)\\
=& \ \left( (\omega_{\nu'}\cdot \omega_\mu)(-1)\right)\cdot \Lambda(1-s,  \widehat{\mathrm g^\circ_{\nu', \nu_n\cdot\abs{\cdot}^s_\sk}(f', \phi_1)} , \widehat f_\mu,  \widetilde{\phi_2}).
\end{aligned}
\end{equation}
Since $s\in \Omega_{\nu, \nu'}$,  Proposition \ref{prpconv1} (b) implies that the three double integrals in \eqref{5.3eq3} are all absolutely convergent.

In view of Proposition \ref{gssur1}, the proposition follows from  \eqref{5.3eq1+2},  \eqref{5.3eq3} and the following lemma.
\end{proof}

\begin{lemp} Assume that $n'=n-1$. Then it holds that % function $\Gamma(s;\nu;\nu')$ defined by \eqref{Gamma} satisfies that
\be \label{5.3eq4}
\Gamma_\psi (s;\nu;\nu')= \Gamma_\psi (1-s; \widehat{\nu'} ;  \widehat{\mu})\cdot  \prod_{1\leq j\leq n-1}  \nu'_j(-1)^n \cdot \prod_{1\leq i, j \leq n-1} \gamma(s, \nu_i\cdot \nu'_j, \psi),
\ee
where $\mu:=(\nu_1,\ldots, \nu_{n-1})$.
\end{lemp}

\begin{proof}
We prove the lemma by induction on $n$. The lemma is easily checked when $n=2$. Assume that $n\geq 3$ and the lemma  holds for $n-1$. Then by the induction hypothesis, we have that
\be \label{5.3eq5}
\begin{aligned}
 \Gamma_\psi (1-s; \widehat{\nu'} ;  \widehat{\mu}) & = \Gamma_\psi (1-s; (\nu'^{-1}_{n-1},\ldots, \nu'^{-1}_1); (\nu^{-1}_{n-1}, \ldots, \nu_2^{-1}))\\
 & =  \Gamma_\psi (s; (\nu_2, \cdots, \nu_{n-1}); (\nu_2',\cdots \nu_{n-2}')) \\
 &\quad  \cdot \prod_{2\leq i\leq n-1} \nu_i(-1)^{n-1}  \cdot  \prod_{2\leq i, j \leq n-1} \gamma(1-s, \nu_i^{-1}\cdot \nu'^{-1}_j, \psi).
 \end{aligned}
\ee

For $\omega\in \widehat{\sk^\times}$, it holds that
\be \label{5.3eq6}
 \gamma(s, \omega,\psi) \cdot\gamma(1-s, \omega^{-1}, \psi) = \varepsilon(s, \omega,\psi) \cdot\varepsilon(1-s, \omega^{-1}, \psi) = \omega(-1).
\ee

By \eqref{5.3eq5} and \eqref{5.3eq6},  the right hand side of \eqref{5.3eq4} equals 
\[
\begin{aligned}
& \quad  \Gamma_\psi (s; (\nu_2, \cdots, \nu_{n-1}); (\nu_2',\cdots \nu_{n-2}'))  \cdot  \prod_{2\leq i\leq n-1} \nu_i(-1)^{n-1} \cdot  \prod_{1\leq j\leq n-1}  \nu'_j(-1)^n  \\
& \cdot \prod_{2\leq i, j \leq n-1}(\nu_i\nu_j')(-1) \cdot \prod_{ 1\leq i, j\leq n-1, \ \min(i,j)=1}   \gamma(s, \nu_i\cdot \nu'_j, \psi),\\
= & \quad  \Gamma_\psi (s; (\nu_2, \cdots, \nu_{n-1}); (\nu_2',\cdots \nu_{n-2}'))  \cdot \nu_1'(-1)^n \cdot \prod_{2\leq i\leq n-1}\nu_i(-1) \\
&\cdot \prod_{ 1\leq i, j\leq n-1, \ \min(i,j)=1}   \gamma(s, \nu_i\cdot \nu'_j, \psi), \\
\end{aligned}
\]
which is easily seen to be equal to $\Gamma_\psi (s;\nu;\nu')$. This finishes  the proof of the lemma.
\end{proof}

Finally, Theorem \ref{thmA} follows  from Propositions \ref{prop5.1}--\ref{prop5.3}.

\section*{Acknowledgements}

D. Liu was supported by the Natural Science Foundation of Zhejiang Province (Grant No.
LZ22A010006) and National Natural Science Foundation of China (Grant No. 12171421). F. Su was supported
by National Natural Science Foundation of China (Grant No. 11901466) and the Qinglan Project of Jiangsu
Province. B. Sun was supported by the National Key $\textrm{R}\,\&\,\textrm{D}$ Program of China (Grant No. 2020YFA0712600).
The authors thank the anonymous referees for the careful reading and comments.

\end{document}